# An Alternating KMF Algorithm to Solve the Cauchy Problem for Laplace's Equation

Chakir Tajani
Faculty of sciences, Ibn Tofail University
Kenitra, Morocco

Jaafar Abouchabaka
Faculty of sciences, Ibn Tofail University
Kenitra, Morocco

## ABSTRACT
This work concerns the use of the iterative algorithm (KMF algorithm) proposed by Kozlov, Mazya and Fomin to solve the Cauchy problem for Laplace's equation. This problem consists to recovering the lacking data on some part of the boundary using the over specified conditions on the other part of the boundary. We describe an alternating formulation of the KMF algorithm and its relationship with a classical formulation. The implementation of this algorithm for a regular domain is performed by the finite element method using the software Freefem. The numerical tests developed show the effectiveness of the proposed algorithm since it allows to have more accurate results as well as reducing the number of iterations needed for convergence.

## 1. INTRODUCTION
In this work, we consider a Cauchy problem for the Laplace equation, called data completion problem, which is to complete the missing data on some part of the boundary (that we cannot assess due to the physical difficulties or inaccessibility geometric) using the over determined data on the other part of the boundary.

This type of problem arises in several areas of engineering such as non destructive control, corrosion detection [1], mechanical problem's particularly in the areas of identification of boundaries on domains, determination of initial condition and fault location [2], in tomography or in Geophysics [3], and also in electroencephalography [4].

The ill-posedness of the problem in the sense of Hadamard makes its resolution by direct methods very difficult, and leads to serious questions including the existence, uniqueness and stability of the solution, that are the three properties required to define well-posed problem according to Hadamard [5].     The existence of the solution of this kind of problem is not always guaranteed, but when the conditions on the accessible part of the boundary are compatible then the existence is assured [6]. Thanks to Holmgreen theorem, we know that this problem has at most one problem [7]. Stability is the most delicate problem since a small perturbation of data provides a large difference between the solution obtained by disturbed data and that obtained by undisturbed data [8]. It suffices here to recall the famous example of Hadamard where he showed for a square domain that with perturbed data, the solution is not bounded even if the data problems tend to zero.

In order to solve the inverse problem for the Laplace equation we have proposed several performing methods to overcome of the ill-posed nature of this kind of problem. The last ancient of them is the one based on optimization tools, introduced by Kohn and Vogelius [9].

Other methods were experimented, among them, we mention the method of Quasi-reversibility introduced by Lattés since 1960 [10], which is to replace the inverse problem by a well-posed problem in the sense of Hadamard by introducing some parameter tends to 0. This method has been adopted by several authors to solve a Cauchy problem, especially Klibanov and Santosa [11], and more recently Bourgeois [12], and others    [13, 14]. This method is effective since it solves the problem directly and the results are accurate and robust. However it has some disadvantages as the particular choice of the parameter introduced which can be difficult to achieve in real circumstances and the difficulty of taking into account the physical constraints that may be related to the problem considered. The regularization method of Tikhonov [15, 16, 17], consists of a minimization problem which is added a penalty term that depends on a parameter called regularization parameter. Tikhonov methods have the disadvantage of disrupting the operator. In addition, these methods need a priori information on the solution of the inverse problem. Other group consists of iterative methods.

Other methods exist like Backus-Gilbert method applied to moment problem [18] and the method applied to the minimization of energy like functional [19].

The group of iterative method has the advantage to allow any physical constraint to be easily taken into account directly in the scheme of the iterative algorithm, simplicity of the implementation schemes and the similarity of schemes for problems with linear and non linear operators. One possible disadvantage of this kind of method is the large number of iterations that may be required in order to achieve convergence.

 Based on these reasons, we have decided in this work to consider the KMF algorithm addressed by Kozlov, Mazya and Fomin  since 1991 [20] (see also [21, 22, 23])., also called alternating method since it is to solve alternately two well-posed problems forming a sequence that approximate the missing conditions problem. To deal with the large number of iterations required to achieve the convergence, relaxation methods have been developed by introducing a relaxation parameter in the Dirichlet condition or the Neumann condition obtained after solving the two well-posed problems of the algorithm [24].     In addition, studies were made on the optimum choice of relaxation parameter to accelerate the convergence of the algorithm [25].

As first work, in order to answer to the question: what data can be introduced to the algorithm to relax and win in number of iterations? we studied in more detail this method, particularly; the relationship between the rate of convergence of this algorithm, the accuracy of the solution and the data of the problem, specially the measure of the inaccessible part of the boundary and the different choice of condition on the accessible part [26]. In the present work, we propose a KMF developed algorithm to reduce the number of iterations needed to achieve convergence with more precision.

The second section is devoted to the presentation of the Cauchy problem for Laplace's equation. Section 3 presents a classical KMF algorithm which enables one to find an approximate solution to that problem. Section 4 considers an alternating KMF algorithm and exhibits the relationship it has with the classical KMF algorithm. Finally, section 5 presents a numerical example showing the feasibility of the alternating formulation, its ability to find an approximate solution more accurately in less iteration.

## 2. MATHEMATICAL FORMULATION
### 2.1 Direct problem statement
Let $\Omega$ be an open set in $\mathbb{R}^2$, with a smooth boundary $\Gamma$. We consider a partition of this boundary $\Gamma = \Gamma_0 \cup \Gamma_1$ where $\Gamma_0 \cap \Gamma_1 = \emptyset$ and $mes(\Gamma_1) \neq 0$.

The direct problem consists to find the harmonic function $u$ solution of the well-posed problem defined as follows:

$$\begin{cases} \Delta u = 0 & on \quad \Omega \\ u = f & in \quad \Gamma_0 \\ \partial_n u = h & in \quad \Gamma_1 \end{cases} \quad (1)$$

Where $f$ and $h$ are the natural boundary conditions.

$\partial_n u$ is the normal derivative of $u$.

This direct problem is well-posed and can be solved by direct method.

### 2.2 Reconstruction inverse problem
In the reconstruction inverse problem, the geometry of the problem is determined, but the boundary conditions are not completely known. This problem arises in cases where a part of the boundary is exposed to environmental conditions which cannot be assessed due to physical difficulties or geometrical inaccessibility. The aim in the reconstruction inverse problem is to find the unknown boundary conditions in $\Gamma_1$ based on the supplementary data provided on the other part $\Gamma_0$.

The problem is to reconstruct a harmonic function u solution of the following problem:

$$\begin{cases} \Delta u = 0 & on \quad \Omega \\ u = f & in \quad \Gamma_0 \\ \partial_n u = g & in \quad \Gamma_0 \end{cases} \quad (2)$$

where $\partial_n u$ is the normal derivative of $u$.

We can notice that no boundary condition is prescribed on the boundary part $\Gamma_1$, and we have two conditions (Dirichlet and Neumann) in the remaining part $\Gamma_0$.

For $f \in H^{\frac{1}{2}}(\Gamma_0)$ and $g \in H^{-\frac{1}{2}}(\Gamma_0)$, where $H^{-\frac{1}{2}}(\Gamma_0)$ is defined as dual space of $H_{00}^{\frac{1}{2}}(\Gamma_0)$ and

$H_{00}^{\frac{1}{2}}(\Gamma_0) = \{v \in L^2(\Gamma_0), \exists w \in H^1(\Omega), w/_{\Gamma_0} = v, w/_{\Gamma_1} = 0\}$

the problem (2) has a unique solution.

## 3. DESCRIPTION OF THE STANDARD ALGORITHM

Consider the Cauchy problem (2) with $f \in H^{\frac{1}{2}}(\Gamma_0)$ and $g \in H^{-\frac{1}{2}}(\Gamma_0)$. The iterative algorithm proposed by Kozlov, Mazya and Fomin (Algorithm standard) investigated is based on reducing this ill-posed problem to a sequence of mixed well-posed boundary value problems and consists of the following steps:

***Step 1.*** Specify an initial guess $u_0 \in H^{\frac{1}{2}}(\Gamma_0)$

***Step 2.*** Solve the following mixed well-posed boundary value problem:

$$\begin{cases} \Delta u^{(0)} = 0 & on \quad \Omega \\ u^{(0)} = u_0 & in \quad \Gamma_1 \\ \partial_n u^{(0)} = g & in \quad \Gamma_0 \end{cases} \quad (3)$$

To obtain $v_0 = \partial_n u^{(0)}/\Gamma_1$ (4)

***Step 3.*** For $n \geq 1$, solve alternatively the two mixed well-posed boundary value problems:

$$\begin{cases} \Delta u^{(2n-1)} = 0 & on \quad \Omega \\ u^{(2n-1)} = v_{n-1} & in \quad \Gamma_1 \\ u^{(2n-1)} = f & in \quad \Gamma_0 \end{cases} \quad (5)$$

To obtain $u_n = u^{(2n-1)}/\Gamma_1$ (6)

and

$$\begin{cases} \Delta u^{(2n)} = 0 & on \quad \Omega \\ u^{(2n)} = u_n & in \quad \Gamma_1 \\ \partial_n u^{(2n)} = g & in \quad \Gamma_0 \end{cases} \quad (7)$$

(8)

To obtain $v_n = \partial_n u^{(2n)}/\Gamma_1$

***Step 4.*** Repeat step.3 from $n > 1$ until a prescribed stopping criterion is satisfied.

### 3.1 Observations
Convergence results were shown by different authors. We cite the results of Kozlov[20], Baumester [27] where it was shown that:

(i) If the Cauchy problem (2) is consistent for the data $(f, g)$; i.e. it has a unique solution $u \in H^1(\Omega)$, then the sequence $(u_n)_{n \geq 0}$ defined in Eq. (6) converges to $u/_{\Gamma_1}$ in the norm of $H^{\frac{1}{2}}(\Gamma_1)$.

(ii) If the sequence $(u_n)_{n \geq 0}$ defined in Eq. (6) converges in $H^{\frac{1}{2}}(\Gamma_1)$ then it converges to $u/_{\Gamma_1}$ where $u \in H^1(\Omega)$ is the solution of the Cauchy problem (2) which in this case exists and is unique.

Furthermore, Jourhmane and Nachaoui [28] showed that:

(iii) If the sequence $(u_n)_{n\geq 0}$ defined in Eq. (6) converges in $H^{\frac{1}{2}}(\Gamma_1)$.then the sequence $(u^{(n)})_{n\geq 0}$ converges in $H^1(\Omega)$ to the solution $u \in H^1(\Omega)$ of the Cauchy problem (2). Moreover, we can notice that:

(iv) The same conclusion is obtained if at the step 1, we specify an initial guess $v_0 \in H^{-\frac{1}{2}}(\Gamma_1)$ instead of an initial guess for $u_0 \in H^{\frac{1}{2}}(\Gamma_1)$, and we modify accordingly the steps 2-3 of the algorithm.

## 4. THE ALTERNATING KMF ALGORITHM

### 4.1 Description of the alternating algorithm

The numerical results obtained by studying the influence of data problems, particularly; the relationship between the measure of the inaccessible part of the boundary and the rate of convergence have shown that:

- The convergence is always guaranteed and it is very fast if the measure of the part of the boundary to be completed is small.
- The convergence requires much more iterations if the inaccessible part is greater.

In this study we develop a new algorithm called "KMF developed Algorithm" in order to improve the rate of convergence of the iterative algorithm described.

The main idea of the algorithm alternative proposed is based on the use of the previous results and the KMF standard algorithm by completing the missing data in alternative way to the two sub-parts of the inaccessible boundary. The inaccessible part is subdivided in two parts, and the KMF standard algorithm is used to complete the data in the first part, then to complete the data in the second part in an alternative way.

For this, we consider $\Gamma_1 = \Gamma_{1,1} \cup \Gamma_{1,2}$ such that $\Gamma_{1,1} \cap \Gamma_{1,2} = \emptyset$ and $mes(\Gamma_{1,1}) = mes(\Gamma_{1,2})$.

The algorithm consists of the following steps:

***Step 1.*** Specify an initial guess $u_0 \in H^{\frac{1}{2}}(\Gamma_1)$

***Step 2.*** Solve the well-posed problem

$$\begin{cases} \Delta u^{(0)} = 0 & on \quad \Omega \\ u^{(0)} = u_0 & in \quad \Gamma_1 = \Gamma_{1,1} \cup \Gamma_{1,2} \\ \partial_n u^{(0)} = g & in \quad \Gamma_0 \end{cases} \quad (9)$$

To obtain $v_{1,0} = \partial_n u^{(0)}_{/\Gamma_{1,1}}$ and $v_{2,0} = \partial_n u^{(0)}_{/\Gamma_{1,2}}$ (10)

***Step 3.*** For $n \geq 1$ solve the two well-posed problems

$$\begin{cases} \Delta u^{(2n-1)} = 0 & on \quad \Omega \\ \partial_n u^{(2n-1)} = v_{1,n-1} & in \quad \Gamma_{1,1} \\ \partial_n u^{(2n-1)} = v_{2,n-1} & in \quad \Gamma_{1,2} \\ u^{(2n-1)} = f & in \quad \Gamma_0 \end{cases} \quad (11)$$

To obtain $u_{1,n} = u^{(2n-1)}_{/\Gamma_{1,1}}$ (12)

And

$$\begin{cases} \Delta v^{(2n-1)} = 0 & on \quad \Omega \\ v^{(2n-1)} = u_{1,n} & in \quad \Gamma_{1,1} \\ \partial_n v^{(2n-1)} = v_{2,n-1} & in \quad \Gamma_{1,2} \\ v^{(2n-1)} = f & in \quad \Gamma_0 \end{cases} \quad (13)$$

To obtain $u_{2,n} = v^{(2n-1)}_{/\Gamma_{1,2}}$ (14)

***Step 4.*** Solve the two well-posed problems

$$\begin{cases} \Delta u^{(2n)} = 0 & on \quad \Omega \\ u^{(2n)} = u_{1,n} & in \quad \Gamma_{1,1} \\ u^{(2n)} = u_{2,n} & in \quad \Gamma_{1,2} \\ \partial_n u^{(2n)} = g & in \quad \Gamma_0 \end{cases} \quad (15)$$

To obtain $v_{1,n} = \partial_n u^{(2n)}_{/\Gamma_{1,1}}$ (16)

And

$$\begin{cases} \Delta v^{(2n)} = 0 & on \quad \Omega \\ \partial_n v^{(2n)} = v_{1,n} & in \quad \Gamma_{1,1} \\ v^{(2n)} = u_{2,n} & in \quad \Gamma_{1,2} \\ \partial_n v^{(2n)} = g & in \quad \Gamma_0 \end{cases} \quad (17)$$

To obtain $v_{2,n} = \partial_n v^{(2n)}_{/\Gamma_{1,2}}$ (18)

***Step 5.*** Repeat step 3. and *step 4.* from $n \geq 1$ until a prescribed stopping criterion is satisfied.

### 4.2 Remarks

If we consider every iteration to consist of solving the four mixed well-posed problems from the Step 3 and 4 of the algorithm, then for every $n \geq 1$ the following approximations are obtained at the iteration number n:

$u_{1,n}$ for the Dirichlet condition on the $\Gamma_{1,1}$.

$u_{2,n}$ for the Dirichlet condition on the $\Gamma_{1,2}$.

$v_{1,n}$ for the Neumann condition on the boundary $\Gamma_{1,1}$.

$v_{2,n}$ for the Neumann condition on the boundary $\Gamma_{1,2}$.

It should be noted that:

- The KMF developed algorithm can be seen as two parallel problems of KMF standard algorithm. These two problems are initialized with the same initial data. Each problem allows to obtain approximation in each subpart $\Gamma_{1,i}$ ou $i = 1,2$ (for the approximation in $\Gamma_{1,1}$ the two well-posed problems (11) and (15), for the approximation in $\Gamma_{1,2}$ the two well-posed problems (13) and (17)).

- Each solved problem allows an approximation in one of the inaccessible subparts that can be introduced in the other well-posed problems.

- The missing Dirichlet condition in the part $\Gamma_1$ can be obtained from the problem (13) since the condition $u_{1,n} = u^{(2n-1)}_{/\Gamma_{1,1}}$ obtained from (11) is introduced in (13) which also provides $u_{2,n} = v^{(2n-1)}_{/\Gamma_{1,2}}$.

➢ The missing Neumann condition in the part $\Gamma_1$ can be obtained from the problem (17), since the condition $v_{1,n} = \partial_n u^{(2n)}_{/\Gamma_{1,1}}$ obtained from (15) is introduced in (17) which also provides $v_{2,n} = \partial_n v^{(2n)}_{/\Gamma_{1,2}}$.

## 5. NUMERICAL RESULTS

In many practical applications $\Gamma_0$ and $\Gamma_1$ are two simple arcs having in common only the end-points.

In this section, we illustrate the numerical results obtained using the alternating KMF algorithm described in section 4. In addition, we investigate the convergence and the accuracy of the solution with respect the number of iterations.

In order to present the performance of the numerical method proposed, we solve the Cauchy problem for an example in a two-dimensional smooth geometry, namely the unit disc $\Omega = \{(x,y) \in \mathbb{R}^2 / x^2 + y^2 \leq 1\}$, also taken in [29] since the condition of a smooth domain is required by the theoretical analysis of Kozlov et al.

We assume that the boundary $\Gamma$ of the solution domain is divided into two disjointed parts $\Gamma_0$ and $\Gamma_1$, namely;

$\Gamma_0 = \{(x,y) \in \mathbb{R}^2 / x = \cos(t), y = \sin(t), \theta \leq t \leq 2\pi\}$
$\Gamma_1 = \{(x,y) \in \mathbb{R}^2 / x = \cos(t), y = \sin(t), 0 \leq t \leq \theta\}$

Where $\theta$ is a specified angle in the interval $(0, 2\pi)$.

In order to illustrate the typical numerical results, we have taken different choices of the angle $\theta$.

The analytical function to be retrieved is given by:

$u(x,y) = x^2 - y^2$

For the implementation of the iterative algorithm we use the software FreeFem, and solve the well-posed problems in the algorithm by the finite element method. In this example, we use a finite element method with continuous piecewise linear polynomials to provide simultaneously the unspecified boundaries Dirichlet and Neumann.

The following stopping criterion was adopted

$E = \|u_n - u_{n+1}\|_{0,\Gamma_1} \leq 10^{-5}$

The convergence of the algorithm may be investigated by evaluating at every iteration the error:

$e_u = \|u_n - u_{ex}\|_{0,\Gamma_1}$ and $e_v = \|v_n - \partial_n u_{ex}\|_{0,\Gamma_1}$

Where $u_n$ is the approximation obtained for the function on the boundary $\Gamma_1$ after n iterations and $u_{ex}$ is the exact solution of the problem of the problem (1). However, in practical applications the error $e_u$ cannot be evaluated since the analytical solution is not known and therefore the error $E$ has to be used.

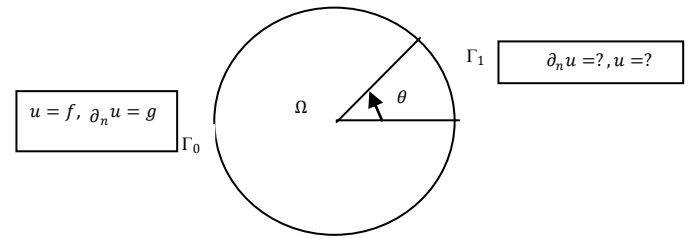

**Fig 1: Schematic diagram showing the prescription of the boundary conditions**

The unknown data in the underspecified boundary in $\Gamma_1$ are given by:

$\partial_n u(x,y) = 2(2x^2 - 1)$ and $u(x,y) = (2x^2 - 1)$

As an initial guess $u_0$ for the step 1 of the algorithm, we have chosen $u_0 = x^2 - x - \frac{1}{2}$.

We notice that $u_0$ is not to close to the analytical solution $u$ on the under specified boundary.

The Fig 2. (resp. fig 3 .) present a comparison between the numerical results $e_u$ ($resp.\ e_v$) obtained with the KMF standard algorithm and the KMF developed algorithm.

It can be seen that the algorithm proposed decreases considerably the number of iteration necessary to achieve the convergence that can be reduced, and present a more accurate approximations for both Dirichlet and Neumann missing data.

For the KMF standard algorithm, the error $e_u$ obtained for $\theta = \frac{\pi}{6}$ after the final iteration 50 is $9.7\ 10^{-3}$. However, with the algorithm proposed it is reduced to $5,6\ 10^{-3}$ after 30 iterations. For $\theta = \frac{\pi}{2}$, the error $e_u$ obtained after 314 iterations is $2.1\ 10^{-2}$. However; it is reduced to $8.\ 10^{-3}$ after 200 iterations.

For the KMF standard algorithm, the error $e_v$ obtained for $\theta = \frac{\pi}{3}$ after the iteration 125 is $5.9\ 10^{-2}$. However, with the KMF developed algorithm, this error is obtained after 68 iterations. In addition, at iteration 125 the error is reduced to $3.6\ 10^{-2}$. For all the results, we can see that the new algorithm allows completing the data more accurately by reducing the number of iterations by two.

Then the new algorithm is more accurate to approximate the Dirichlet and Neumann conditions in the inaccessible part.

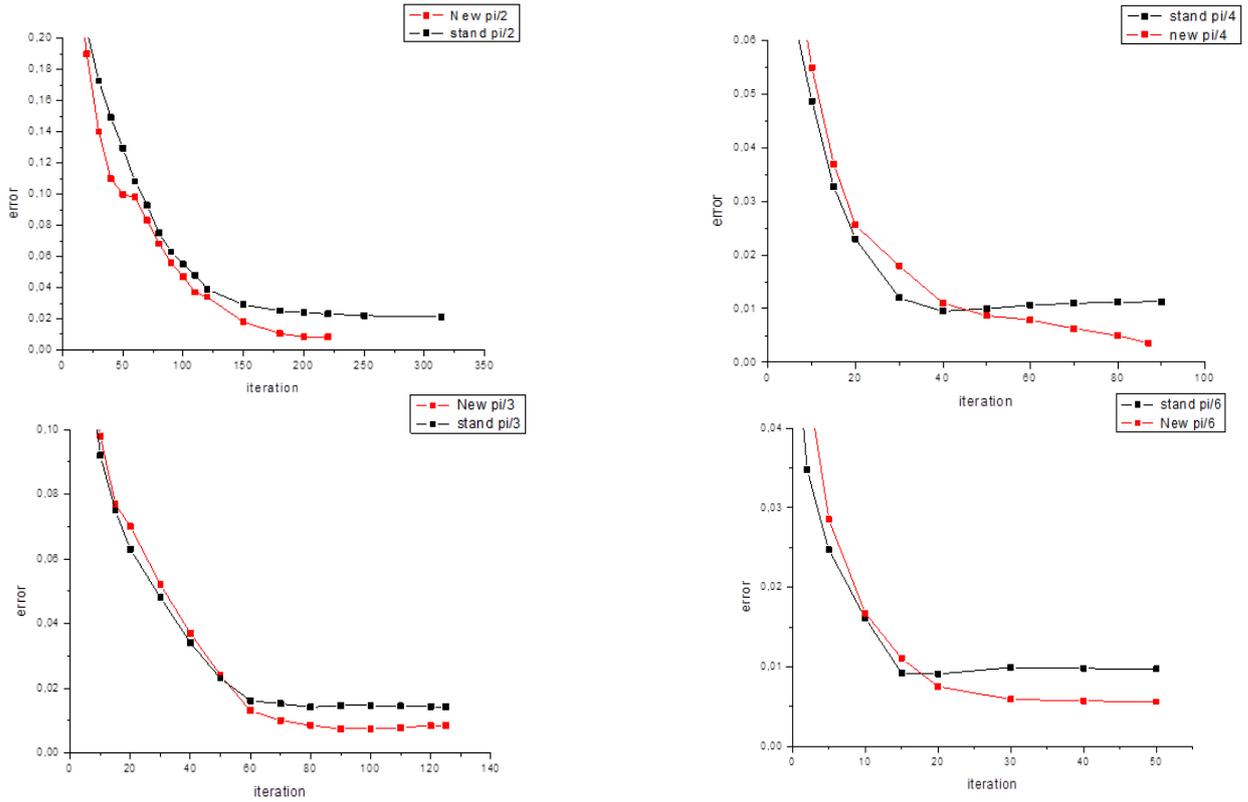

**Fig 2**: The error $e_u$ as a function of the number of iterations for different choice of $\theta$ ($\frac{\pi}{2}, \frac{\pi}{3}, \frac{\pi}{4}, \frac{\pi}{6}$) obtained for KMF developed algorithm (New) in comparison with classical algorithm (Stand)

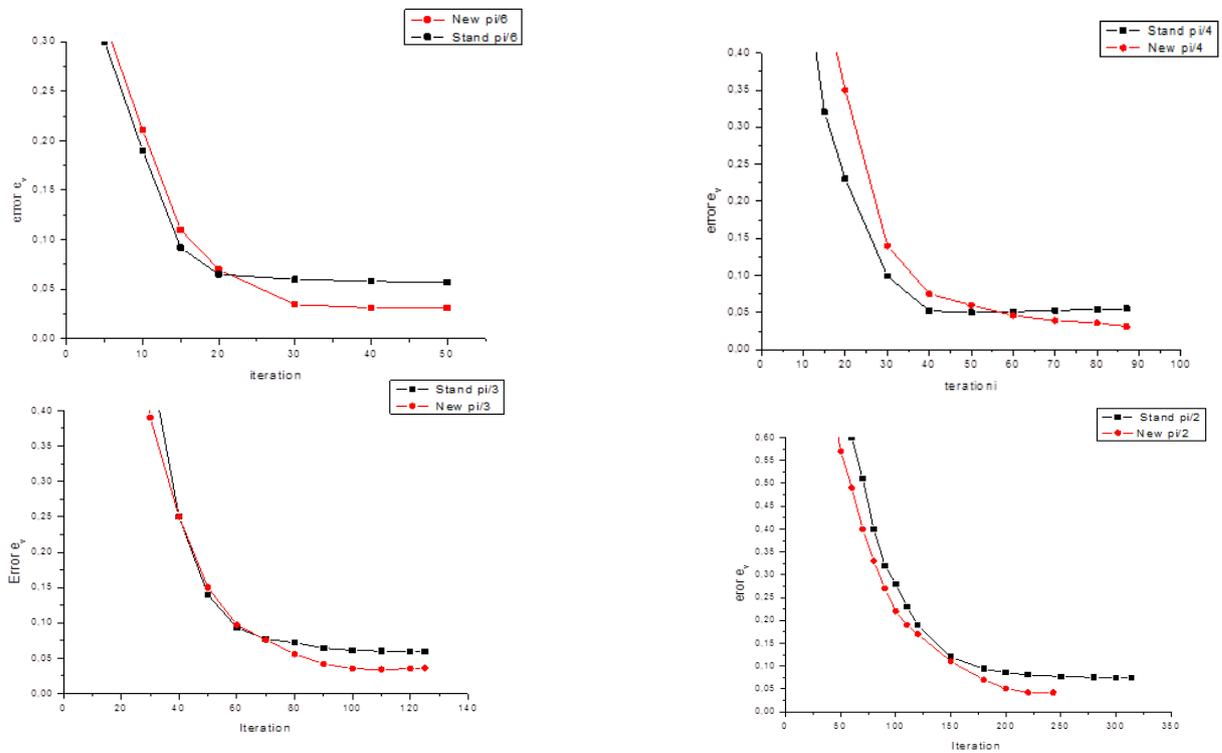

**Fig 3**: The error $e_v$ as a function of the number of iterations for different choice of $\theta$ ($\frac{\pi}{2}, \frac{\pi}{3}, \frac{\pi}{4}, \frac{\pi}{6}$) obtained for KMF developed algorithm (New) in comparison with classical algorithm (Stand)

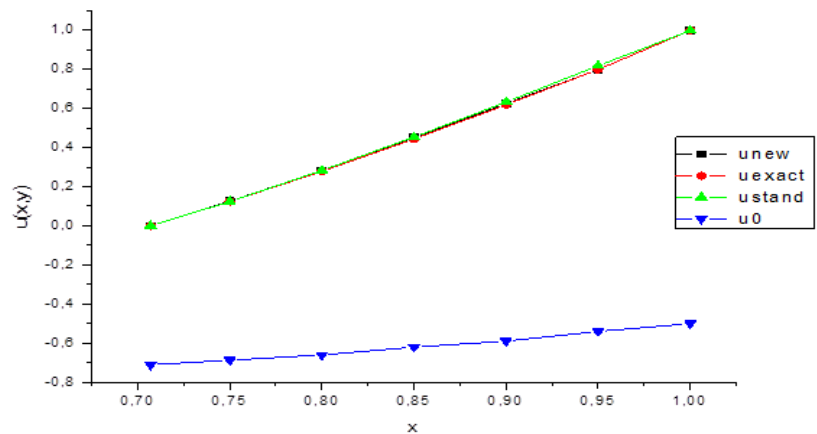

**Fig 4: The numerical results for the function u on the boundary Γ$_1$ obtained with KMF developed algorithm (unew) in comparison with the analytical solution (uexact), the initial guess (u$_0$) and the solution with KMF standard algorithm (ustand)**

Fig 4. shows the numerical results obtained in approximating the function u in the part of the boundary Γ$_1$, indicating that from a choice of an initial data, we obtain satisfying results for both algorithms. However, the KMF developed algorithm requires less iteration to achieve more accurate convergence.

Such implementation of the algorithm allows us to notice that after the number of iterations is sufficiently increased, the error become small; this shows that the numerical solution is accurate and consistent with the number of iterations.

Furthermore, when perturbations are introduced into the given data problem the numerical results obtained are stable. Overall, it can be concluded that the alternating iterative algorithm (called KMF developed algorithm) proposed produces a convergent, stable and accurate numerical solution.

## 6. CONCLUSION

In this paper we have investigated the KMF iterative algorithm for a Cauchy problem for Laplace equation. Comparison of numerical results with those obtained by the KMF standard algorithm show that the proposed algorithm significantly reduces the number of iterations needed to achieve the convergence and produces more accurate results. In addition, it can be concluded that the proposed algorithm is very efficient to reduce the rate of convergence. When perturbations are introduced into the given date problem the results are stable. Overall, it can be concluded that the alternating iterative algorithm proposed produces a convergent, stable and accurate numerical solution.

In this paper, the part of the boundary is divided in two parts. Work is in progress for implementing this algorithm by dividing the inaccessible part in more than two parts to more accelerate the convergence and implementing the new algorithm in the case of the Cauchy problem for Helmholtz equation.